\documentclass[10pt]{article}

\renewcommand{\Box}{\square}

\newcommand{\isd}{\stackrel{d}{=}}

\newcommand{\RR}{\mathbb{R}}
\newcommand{\PP}{\mathbb{P}}
\newcommand{\ZZ}{\mathbb{Z}}
\newcommand{\NN}{\mathbb{N}}
\newcommand{\EE}{\mathbb{E}\,}

\newcommand{\Ber}{\textrm{Ber}}

\newcommand{\tS}{{\tilde{S}}}
\newcommand{\tD}{{\tilde{D}}}

\newcommand{\uA}{{\underline{A}}}
\newcommand{\uD}{{\underline{D}}}

\usepackage{amsfonts}
\usepackage{amsxtra}
\usepackage{amsmath}
\usepackage{a4}
\usepackage{amssymb}
\usepackage{theorem}
\usepackage{psfrag}
\usepackage{epsfig}
\usepackage{graphicx}
\usepackage{color}

\newtheorem{theorem}{Theorem}
\newtheorem{lemma}[theorem]{Lemma}
\newtheorem{corollary}[theorem]{Corollary}
{\theorembodyfont{\rmfamily}}
{\theorembodyfont{\rmfamily}}
\newtheorem{proposition}[theorem]{Proposition}
\numberwithin{theorem}{section}
\numberwithin{equation}{section}
\numberwithin{figure}{section}

\renewcommand{\Ber}{\ensuremath\textrm{Ber}}
\newcommand{\Geom}{\ensuremath\textrm{Geom}}

\allowdisplaybreaks

\begin{document}

\title{Fixed points for multi-class queues}
\author{
\textbf{James B.~Martin and Balaji Prabhakar}
\\
\textit{University of Oxford and Stanford University}}
\date{ }
\maketitle


\begin{abstract}
Burke's theorem can be seen as a fixed-point result for 
an exponential single-server queue; when the arrival
process is Poisson, the departure process has the same distribution
as the arrival process. We consider extensions
of this result to multi-type queues, in which different types of customer
have different levels of priority. We work with a model of a queueing server
which includes discrete-time and continuous-time $M/M/1$ queues
as well as queues with exponential or geometric service batches
occurring in discrete time or at points of a Poisson process. The fixed-point
results are proved using \textit{interchangeability} properties for queues in tandem,
which have previously been established for one-type $M/M/1$ systems. 
Some of the fixed-point results have previously been derived as a consequence of
the construction of stationary distributions for multi-type interacting particle systems,
and we explain the links between the two frameworks. The fixed points have interesting
``clustering'' properties for lower-priority customers. An extreme case is an example of a Brownian 
queue, in which lower-priority work only occurs at a set of times of measure 0 
(and corresponds to a local time process for the queue-length process of higher priority work). 
\end{abstract}

\section{Introduction}
One of the most famous results in queueing theory is \textit{Burke's theorem} \cite{Burke}.
Consider a queue in which available services occur as a Poisson process of rate $\mu$
(a so-called $./M(\mu)/1$ queueing server).
If the arrival process is a Poisson process of rate $\lambda<\mu$ 
(independent of the service process), then the departure process
is also a Poisson process of rate $\lambda$. We may say 
that the arrival process is a \textit{fixed point} for the server.

In this paper we consider the question of fixed points for queues 
with two or more classes of customer (with different levels of priority). 
When a service occurs in such a queue, it is used by a customer
whose priority is highest out of those currently present in the queue. 
We will see that a two-type fixed point can be constructed using the output processes
(consisting of departures and unused services) from a one-type queue.
Then in a recursive way, a fixed point with $m\geq3$ classes 
of customer can be constructed using the output of a queue whose arrival
process is itself a fixed point with $m-1$ classes.

Except in the familiar one-type case, the fixed points are not Markovian.
In particular, one observes \textit{clustering} of the lower-priority customers.

In the paper we work with a queueing model which is somewhat more general
than the $./M/1$ queue described above. Our basic model is 
of a discrete-time queue with batch arrivals and services.
Let $S_n$ be the amount of service offered at time $n$. 
We obtain fixed-point results for the case where $S_n$ are i.i.d.\ 
and each $S_n$ has so-called ``Bernoulli-geometric'' distribution,
i.e.\ is equal to the product of a geometric random variable
and an independent Bernoulli random variable. By taking appropriate
limits where necessary, this model covers a variety of
previously considered queueing servers, for example discrete-time $./M/1$ 
queues \cite{HsuBurke}, discrete-time queues with geometric
or exponential service batches \cite{BedAzi, DraMaiOco, OConnell},
continuous-time $./M/1$ queues as described above,
continuous-time queues with geometric or exponential service batches
occurring at times of a Poisson process, and Brownian queues \cite{Harrisonbook,
OcoYor}. Versions of Burke's theorem and related reversibility results
were proved for this Bernoulli-geometric model in \cite{jbmbatch}.

Some such fixed-point processes were already constructed in certain 
cases
($M/M/1$ queues in continuous or discrete time) in 
\cite{femaihp} in the context of stationary distributions for certain 
multi-type interacting particle
systems. In this paper we give a more direct proof of the 
fixed-point property, which relies on properties 
of \textit{interchangeability} for queueing servers. Weber \cite{Weberinterchangeability}
showed that for a tandem queueing system consisting of 
two independent $./M/1$ servers with service rates $\mu_1$ and $\mu_2$,
and an arbitrary arrival process, the distribution of the departure
process is unchanged if $\mu_1$ and $\mu_2$ are exchanged. 

This interchangeability 
result was subsequently proved in a number of different ways, for example 
in \cite{Venkatinterchangeability}, \cite{Lehtoneninterchangeability}, and 
\cite{TsoWalinterchangeability}. The coupling proof given by 
Tsoucas and Walrand in \cite{TsoWalinterchangeability} is 
important for our purposes, since we can use their approach to extend the 
interchangeability
result to multi-type queues.

Before developing the general batch queueing model, we begin in Section 
\ref{dotM1section} 
by giving a guide to the main results and methods of proof in 
the particular case of the continuous-time $./M/1$ queue. 
Since this model is already rather
familiar, we give an informal account without
introducing too much notation. (Everything is developed rigorously 
in later sections). In addition, certain aspects are simpler in the
$./M/1$ case; for example, the service process has only one parameter,
so all such service processes are interchangeable, and given any 
vector of arrival rates $\lambda_1,\dots,\lambda_m$ (corresponding
to customers of classes $1,\dots,m$ respectively), there is a unique 
fixed-point arrival process which is common to all $./M(\mu)/1$ queues
for $\mu>\lambda_1+\dots+\lambda_m$. 

Our general model is introduced in Section \ref{model}, which 
describes the set-up of a discrete-time
batch queue. Multi-class systems are introduced in Section 
\ref{multiclass}. The Bernoulli-geometric distribution,
and corresponding queueing servers, are described in Section
\ref{BerGeom}.

Interchangeability results are given in Section \ref{interchangeability}.
These extend the results for one-type $./M/1$ queues described above,
to cover multi-type systems and to the more general queueing server model.
In Section \ref{mainresultsection} we give the 
construction of multi-type fixed points, and prove the fixed-point
property using the interchangeability results. The main 
result is given in Theorem \ref{fixedpointtheorem} (the corresponding
results in the $./M/1$ case are 
Theorem \ref{2typeMM1thm} and Theorem \ref{mtypeMM1thm}). 
The proof
of the interchangeability result itself is given in Section \ref{interchangeabilityproof}.

In Section \ref{examplesection} we give examples of the application of the results to 
several of the particular queueing systems described above. The final example
is that of the Brownian queue. Here the lower-priority work in the fixed-point process
corresponds to the local-time process of a reflecting Brownian motion; 
this process is non-decreasing and continuous but is constant except on a 
set of measure 0. This is an extreme case of the ``clustering of lower-priority customers''
referred to above.

The connections with interacting particle systems are discussed in 
Section \ref{particlesection}. The fixed points for $M/M/1$ servers
in discrete time and in continuous time correspond to stationary
distributions for multi-type versions of the TASEP \cite{FerMarmulti}
and of Hammersley's process \cite{FerMarHAD}, respectively. Time in the queueing
systems corresponds to space in the particle systems; with this identification,
questions of fixed points for queues and stationary distributions for particle
systems are closely analogous. 

Finally in Section \ref{continuoussection} we mention a limit as the number of classes goes 
to infinity, with the density of each class going to 0. In this limit, the class-label of 
each customer becomes, for example, a real number in $[0,1]$. This gives another illustration 
of the clustering phenomenon; although a priori, any given label has probability 0 of 
occurring, nonetheless for any realisation of the process, each label that does occur will 
occur infinitely often with probability 1.

\section{Continuous-time $./M/1$ case} 
\label{dotM1section}
\subsection{$./M/1$ queueing servers and Burke's theorem}

Let $A$ and $S$ be independent Poisson processes, of rate $\lambda$ and $\mu$ respectively,
with $\lambda<\mu$. 
We can use these processes to define an $M/M/1$ queue with arrival rate $\lambda$ and
service rate $\mu$. 

Arrivals occur at points of the process $A$; at these points the queue-length increases by 1.
At a point of the process $S$, a departure occurs and the queue-length decreases by 1,
unless the queue-length is already 0, in which case it stays the same and
we say that an \textit{unused service} has occurred.

We write $D$ for the process of departures and $U$ for the process of unused services.

Perhaps the most famous result in queueing theory is \textit{Burke's theorem},
which states that the departure process $D$ is itself a Poisson process of rate $\lambda$.
We may regard the service process $S$ as an operator (a ``$./M/1$ queue", or a 
``$./M(\mu)/1$ queue'' if we want to emphasise the service rate) which maps the distribution 
of an arrival process to the distribution of a service process. In this sense, 
Burke's theorem \cite{Burke} says that 
a Poisson process of rate $\lambda$ is a \textit{fixed point} for a $./M(\mu)/1$ queue
whenever $\mu>\lambda$. 

In \cite{Anantharamuniqueness}, Anantharam showed that in fact Poisson processes
are the only such ergodic fixed points. Mountford and Prabhakar \cite{MouPra} showed further 
that these fixed points are \textit{attractive}; if one starts with any 
ergodic process of rate $\lambda$ and applies the ``$./M(\mu)/1$ queueing operator''
repeatedly, the sequence of distributions obtained converges weakly to a Poisson process
of rate $\lambda$.

\subsection{Multi-class queues}

In this paper we consider fixed-points in the context of \textit{multi-class} queues. 

Consider again a queue whose service process is a Poisson process of rate $\mu$. 
The queue may now contain several types of customer, say types $1,2,\dots,m$.
Each arrival to the system is of one of these types.
Customers of type 1 (``first-class customers'') have the highest priority,
followed by those of type 2 (``second-class customers'') and so on.
When a service occurs in the queue, if there are any customers present, then
a customer with the highest priority out of those present departs from the system. 
Hence each departure from the system also has a type. Again one may have 
``unused services'', when the queue is completely empty at the time of an event
in the service process.

An $m$-type queue may be seen as a coupling of $m$ one-type queues,
which share the same service process. Given $r$ with $1\leq r\leq m$, 
consider only customers of classes 1 up to $r$, ignoring the differences
between these customers. Since any such customer has priority any customer 
whose class is higher than $r$, the process obtained behaves exactly
as a one-type queue. 

\subsection{Multi-class fixed points}\label{MM1multiclass}

The distribution of a multi-type arrival process is said to be a \textit{fixed point}
if the process of departures from the system has the same distribution as the process of 
arrivals. (We assume that arrivals and services are independent).

We will consider arrival processes which are stationary and ergodic, in which case
for each $m$ there is a deterministic long-run intensity of arrivals of customers 
of type $m$. 

A coupling argument analogous to that used by Mountford and Prabakhar in \cite{MouPra} 
can be used to 
show that for any $\lambda_1,\dots,\lambda_m$ with $\lambda_1+\dots+\lambda_m<\mu$, there is a unique 
stationary and ergodic $m$-type arrival process with intensity $\lambda_r$ of 
customers of type $r$ which is a fixed point for the $./M(\mu)/1$ queue. 

As observed above, the system comprising only customers of types 1 up to $r$ 
can be seen as a single one-type queue (for each $1\leq r\leq m$). 
Since the only one-type fixed points are 
Poisson processes, this shows that in any multi-type fixed point, 
the combined process of customers of types $1,\dots,r$ must be a Poisson
process of rate $\lambda_1+\dots+\lambda_r$.

\subsection{Construction of a 2-type fixed point}
First we describe how to construct such a two-type fixed point, 
with intensities $\lambda_1$ and $\lambda_2$.

Consider the output process $(D,U)$ of an $M/M/1$ queue as described above, 
with arrival rate $\lambda_1$ and service rate $\lambda_1+\lambda_2$.
We regard $(D,U)$ as a two-type process, in which first-class customers
occur at the points of $D$ and second-class customers occur at the points of $U$. 
Write $F_2=F_{2,\lambda_1,\lambda_2}$ for the distribution of this process.
\begin{theorem}\label{2typeMM1thm}
The distribution $F_{2,\lambda_1,\lambda_2}$
is a fixed point for a $./M/1$ queue with 
service rate $\mu$, for any $\mu>\lambda_1+\lambda_2$.
\end{theorem}

Note that indeed the sub-process of first-class customers is a Poisson 
process of rate $\lambda_1$ (by Burke's theorem) and that the 
combined process of first-class and second-class customers is a Poisson
process of rate $\lambda_1+\lambda_2$ (since $D+U=S$). 

However, note also that the process $U$ of second-class customers is \textit{not} a Poisson process. 
Rather, the second-class customers tend to \textit{cluster}. 
For example, fix $\lambda_1$ and suppose that $\lambda_2$ is very small. The process $(D,U)$
is obtained as the output of a queue whose arrival rate is $\lambda_1$ and whose service
rate is $\lambda_1+\lambda_2$, so that it is operating very close to capacity. 
The long-run rate of unused services is only $\lambda_2$. However, if we observe
an unused service, we know that the queue is currently empty; hence
just after an unused service, the instantaneous rate at which another unused service
occurs is $\lambda_1+\lambda_2$ which may be much larger. 

\subsection{Recursive construction of multi-type fixed points}
Now we show how to construct 
fixed points with a larger number of classes recursively. 
Fix $\lambda_1,\dots,\lambda_m$. 
Suppose we have already constructed a distribution of an $(m-1)$-type process 
$F_{m-1}$, with intensity $\lambda_i$ of $i$th-class customers for $1\leq i\leq m-1$,
which is a fixed point 
for a $./M(\mu)/1$ queueing server (whenever $\mu>\lambda_1+\dots+\lambda_{m-1}$).
Now consider a $./M(\lambda_1+\dots+\lambda_m)/1$ 
queue whose arrival process has distribution $F_{m-1}$.
Write $(D_1,D_2,\dots,D_{m-1},U)$ for the output process of this queue, 
comprising departures of types $1,2,\dots,m-1$ and unused services. 
Now identify points of $U$ as $m$th-type customers, and write $F_m=F_{m,\lambda_1,\dots,\lambda_m}$
for the distribution of the $m$-type process obtained. 
\begin{theorem}\label{mtypeMM1thm}
The distribution $F_{m,\lambda_1,\dots,\lambda_m}$ is a fixed point for a 
$./M/1$ queue with service rate $\mu$,
whenever $\mu>\lambda_1+\dots+\lambda_m$.
\end{theorem}

\begin{figure}[t]
	\centering
	{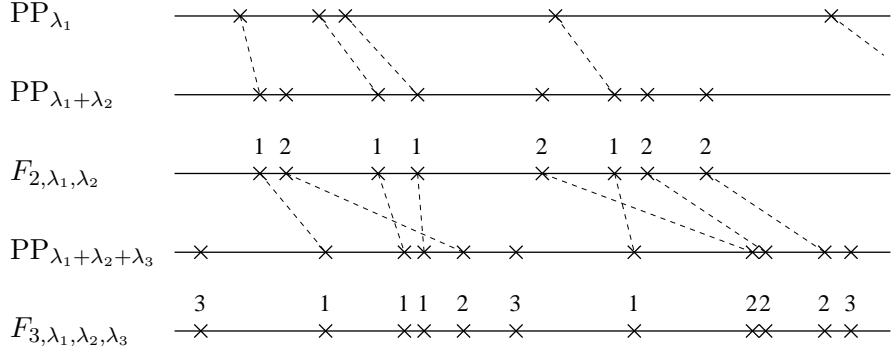}
	\caption{The process $F_2$ is constructed using the departure and unused service
processes of a one-type queue whose arrivals are Poisson rate $\lambda_1$ and whose services 
are Poisson rate $\lambda_1+\lambda_2$. The process $F_3$ is then constructed
using the departure and unused service processes of a two-type queue whose arrivals
have distribution $F_2$ and whose services are Poisson rate $\lambda_1+\lambda_2+\lambda_3$.
\label{linesfig}}
\end{figure}

The construction of $F_2$ and $F_3 $ is illustrated in Figure \ref{linesfig}.

Note that if we take the a process with distribution $F_m$ and ignore the $m$th-class customers, 
we obtain a process with distribution $F_{m-1}$. (This had to be the case, 
because of the uniqueness of the fixed-points with given intensitites, 
and because the ignoring the $m$th class of an $m$-class system simply gives 
an $m-1$-class system). 

Theorem \ref{mtypeMM1thm} was already proved in \cite{FerMarHAD}; 
it emerged as a corollary of a multi-type tandem queue construction 
which was used to give the stationary distribution of 
a multi-type version of an interacting particle system called ``Hammersley's process''.
In this paper we give a much more direct proof, using 
properties of \textit{interchangeability} of queueing servers. 

\subsection{Tandems and interchangeability}
In 1979 Richard Weber \cite{Weberinterchangeability} proved an 
\textit{interchangeability} property for $./M/1$ queueing servers.

Consider two independent $./M/1$ queueing servers in tandem, 
with service rates $\mu_1$ and $\mu_2$. The first queue has some arrival process $A$,
with an arbitrary distribution (for example, it could even be deterministic), 
which is independent of the service processes. By ``tandem'' we mean that 
a customer leaving the first queue immediately joins the second queue;
the departure process from the first server is the arrival process of the second. 

The result of \cite{Weberinterchangeability} is that the law of the departure
process from the system (that is, the departure process from the second queue),
is unchanged if $\mu_1$ and $\mu_2$ are exchanged. 

Given a queue with arrival process $A$ and service process $S$, we
write $D(A,S)$ for the departure process from the queue. So in a system
with two queues in series, with arrival process $A$ and service process $S_1$
and $S_2$ at the first and second queues respectively, the process of departures
from the system is $D(D(A,S_1),S_2)$. Then the result of \cite{Weberinterchangeability}
can be written as follows:
\begin{theorem}\label{Weberthm}
Let $S_1$ and $S_2$ be independent Poisson processes of rates $\mu_1$ and $\mu_2$ 
respectively. Then for any arrival process $A$,
\[
D(D(A,S_1), S_2)\isd D(D(A,S_2), S_1).
\]
\end{theorem}
Thus, for any arrival process, the order of the queues does not affect the 
law of the output of the system. By induction, the result extends easily to tandems
containing more than two queues.

Alternative proofs of this interchangeability result were subsequently given by 
Anantharam \cite{Venkatinterchangeability}, by 
Lehtonen \cite{Lehtoneninterchangeability}, and by Tsoucas and Walrand \cite{TsoWalinterchangeability}. 
The proof by Tsoucas and Walrand in \cite{TsoWalinterchangeability}
is particularly important for our purposes since it allows an extension 
of the result which will also apply to multi-class systems. 
Their methods can be used to give the following result.
\begin{theorem}\label{MM1couplingthm}
Let $S_1$ and $S_2$ be independent Poisson processes of rates $\mu_1$ and $\mu_2$.
There is a coupling of $S_1$ and $S_2$ with two further processes $\tS_1$ and $\tS_2$
such that:
\begin{itemize}\item[(i)]$(S_1,S_2)\isd(\tS_2,\tS_1)$; and
\item[(ii)] If $A$ is any arrival process, then $D(D(A,S_1),S_2)=D(D(A,\tS_1), \tS_2)$. 
\end{itemize}
\end{theorem}

The key point is that we can couple a system with service rates $\mu_1,\mu_2$ 
and another system with service rates $\mu_2,\mu_1$ in such a way 
that \textit{for any arrival process}, the departure processes from the two systems
are the same. (Weber's result already implies the existence of such a coupling 
for any fixed arrival process, but not necessarily that the same coupling works 
for all arrival processes simultaneously).

Above we explained how a queue with $m$ types of customer can be regarded 
as a coupling of $m$ one-type systems, with the same
service process but with different arrival processes. 
(In the $r$th system, we only look at customers whose
class in the original system is between $1$ and $r$, and neglect all the customers
whose class is greater than $r$). Since Theorem \ref{MM1couplingthm} 
gives a coupling which yields interchangeability simultaneously for 
all arrival processes, it immediately implies interchangeability for any multi-type arrival process:
\begin{corollary}\label{MM1multitype}
Let $S_1$ and $S_2$ be independent Poisson processes of rates $\mu_1$ and $\mu_2$.
Then for any multi-type arrival process $\uA$, 
\[
D(D(\uA,S_1), S_2)\isd D(D(\uA,S_2), S_1).
\]
\end{corollary}

The coupling described in Theorem \ref{MM1couplingthm} is quite easy to construct. 
(The construction in \cite{TsoWalinterchangeability} is essentially the same
although not expressed in quite the same way). Assume without loss of generality
that $\mu_1<\mu_2$. Now consider a queue whose arrival process is $S_1$ and whose
service process is $S_2$. Let $D=D(S_1,S_2)$ be the departure process 
of this queue, and let $U$ be the process of unused services (so that
$S_2=D+U$). 

Now let $\tS_2=D$, and let $\tS_1=S_1+U$. Loosely speaking, 
we have transferred the unused service events from one process to the other. 

The fact that $(\tS_2, \tS_1)$ has the same distribution as $(S_1,S_2)$ is
an extension of Burke's Theorem, and can be proved using simple reversibility arguments
(see for example Theorem 3 of \cite{neilsurvey}). 
To prove the interchangeability statement in part (ii) of 
Theorem \ref{MM1couplingthm}, one can verify pathwise that, whatever the arrival process is, 
transferring ``unused services'' in this way cannot affect the overall departure process.

A full proof of this result in the more general setting is given below
(Theorem \ref{interchangethm}). 
The construction in the general setting is slightly more complicated; 
again we remove the process of ``unused services'' from $S_2$ to obtain $\tS_2$,
but now we do not add precisely the same process to $S_1$ to obtain $\tS_1$,
but instead something more like an independent copy of the process. 
However, the same ideas of reversibility are again the key part of the proof. 

\subsection{Proving multi-class fixed points using interchangeability}
Still in the context of the $./M/1$ queue, we finally indicate 
how the interchangeability result of Corollary \ref{MM1multitype}
can be used to prove that the processes $F_m$ defined above
 are indeed fixed points. The arguments
in this section are informal but can be made fully rigorous, as they 
are below in the proof of Theorem \ref{fixedpointtheorem}.
We will show the argument for $m=2$; the case $m>2$ follows
by induction in a natural way, and can be found in the general setting below. 

Recall that $F_2=F_{2,\lambda_1,\lambda_2}$ is the distribution of a two-type process which 
can be obtained by taking the departure and unused service processes $(D,U)$ from 
an $M/M/1$ queue whose arrival rate is $\lambda_1$ and whose service rate is 
$\lambda_2$. So the first-class customers occur as a Poisson process of rate
$\lambda_1$, and the first-class and second-class customers together occur
as a Poisson process of rate $\lambda_1+\lambda_2$.

Let $\mu>\lambda_1+\lambda_2$, and consider 
feeding the arrival process $F_2$ into a system of two queues in series,
with service rates $\lambda_1+\lambda_2$ and $\mu$. The two possible orderings
of the two queues 
can be visualized
as follows:
\begin{align*}
F_2\longrightarrow\boxed{S_{\lambda_1+\lambda_2}}&\longrightarrow \makebox[1em]{$G$} \longrightarrow
\boxed{\hspace{3.23mm}{S_\mu}\hspace{3.23mm}}
\longrightarrow H
\\
F_2\longrightarrow
\boxed{\hspace{3.23mm}{S_\mu}\hspace{3.23mm}}
&\longrightarrow \makebox[1em]{$J$} \longrightarrow
\boxed{S_{\lambda_1+\lambda_2}}
\longrightarrow K
\end{align*}
For example, $G$ is the distribution of the departure process 
of the first queue with service rate $\lambda_1+\lambda_2$,
and is then itself used as the arrival process for the second queue with service rate $\mu$.

Note that the total arrival rate in the process $F_2$ is $\lambda_1+\lambda_2$. 
But the first queue in the first system only has service rate $\lambda_1+\lambda_2$,
so it is saturated; every service is used by some customer. 
Since the arrival process of first-class customers is a Poisson process 
of rate $\lambda_1$, the first-class departures in $G$ occur as the departures
from an $M(\lambda_1)/M(\lambda_1+\lambda_2)/1$ queue; the second-class customers
fill up all the other services which are unused by first-class departures. 
Hence this two-type process 
has the distribution of the departure and unused service process $(D,U)$ from an
$M(\lambda_1)/M(\lambda_1+\lambda_2)/1$ queue, and we have $G=F_2$.

The same argument shows also that $K=F_2$, since the process $J$ has total intensity 
$\lambda_1+\lambda_2$, and (by Burke's theorem) its first-class customers 
occur as a Poisson process of rate $\lambda_1$.

Finally, the interchangeability result in \ref{MM1multitype} tells us that $H=K$,
so that also $H=F_2$.
Since we have $G=H=F_2$, we can concentrate only 
on the second queue in the first line, 
to obtain
\[
F_2\longrightarrow
\boxed{\hspace{3.23mm}{S_\mu}\hspace{3.23mm}}
\longrightarrow F_2
\]
so that indeed $F_2$ is a fixed point for the $./M(\mu)/1$ queue as desired. 

\section{Batch queueing model}\label{model}
We now move to the more general queueing model. We begin by defining
a model of a batch queue in discrete-time.
By taking particular values, or appropriate limits, this model will cover
various interesting cases, including 
discrete-time and continuous-time $M/M/1$ queues.

The batch queue is driven by an arrival process $\left(A(n), n\in\ZZ\right)$ and a 
service process $\left(S(n), n\in\ZZ\right)$.

At time-slot $n\in\ZZ$, $A(n)$ customers arrive at the queue. 
Then service is available for $S(n)$ customers;
if the queue-length is at least $S(n)$, then $S(n)$ customers are served, while 
if the queue length is less than $S(n)$ then all the customers are served. 

We define various processes as functions 
of the basic data $A$ and $S$ of the queue. 

Let $X(n)$ be the queue length after the service $S(n-1)$, before the arrival $A(n)$.
Formally, we define the process $X$ by 
\begin{equation}\label{Xdef}
X(n)=\sup_{m\leq n}\sum_{r=m}^{n-1}\left(A(r)-S(r)\right)
\end{equation}
(where a sum from $n$ to $n-1$ is understood to be 0).
In accordance with the description of the queue above, 
we have the basic recurrences
\begin{equation}\label{Xrec}
X(n+1)=\left[ X(n) + A(n) - S(n)\right]_+
\end{equation}
(where $[x]_+$ denotes $\max\{x,0\}$).

Let $D(n)$ be the number departing from the queue at the time of the service $S(n)$. So
\begin{align*}
D(n)&=\min(X(n)+A(n), S(n))\\
&=X(n)+A(n)-X(n+1).
\end{align*}

Finally let $U(n)=S(n)-D(n)$ be the unused service at the time of the service $S(n)$.

\begin{figure}[t]
	\centering
	{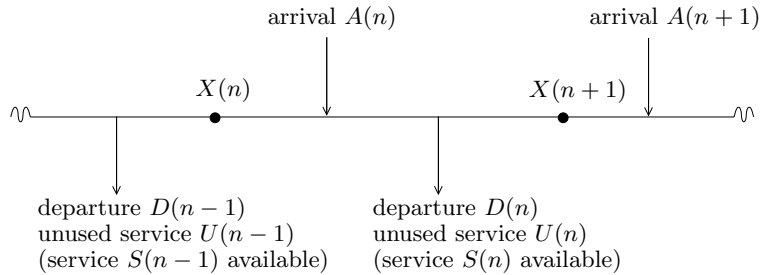}
	\caption{The evolution of the queue with batch services
and arrivals \label{slotfig}}
\end{figure}

See Figure \ref{slotfig} for a representation of the evolution 
of the queue along with its inputs and outputs. 

Since $X$, $D$, $U$ are all functions of the data $A$ and $S$, we sometimes write
$D=D(A,S)$, and so on. 

Note that we allow the possibility that $X(n)=\infty$. 
Indeed, we don't impose any stability condition on the queue; so, for example, 
if the average rate of arrivals exceeds the average
rate of service, then the queue will become saturated. In fact, the following
simple observation will be useful later:
\begin{lemma}\label{saturationlemma}
Suppose the arrival and service processes are independent, 
and $A_n$ are i.i.d.\ with mean $\lambda$ and $S_n$ are i.i.d.\ with mean $\mu\leq\lambda$.
Then with probability 1, 
\begin{align*}
X_n&=\infty \text{ for all }n, \\
D_n&=S_n \text{ for all }n, \\
U_n&=0 \text{ for all }n.
\end{align*}
\end{lemma}
The lemma follows immediately from the definition 
(\ref{Xdef}), since the random walk whose increment at step $r$ is $A(-r)-S(-r)$
is either recurrent (if $\mu=\lambda$) or escapes to $+\infty$ with probability 1 
(if $\mu<\lambda$); in either case it attains arbitrarily high values.  

If the $A(n)$ and $S(n)$ take integer values, then it is natural to talk in terms
of ``number of customers'' arriving, or departing, or in the queue, and so on.
However we will also consider cases where the values are more general, in which case
one could talk of ``amount of work'' rather than ``number of customers''.

\section{Multi-class queues}\label{multiclass}
We will now define a \textit{multi-class} batch queue. The system can now contain
different types of customers (or work) with different priorities. When service occurs
at the queue, it is first available to first-class customers. If there is more service available
than there are first-class customers present in the queue, the remaining service 
(unused by the first-class customers) is then offered to customers of lower class,
starting with second-class customers, then third-class and so on. 

For example, suppose that at the start of time-slot $n$, there are 7 customers in the queue,
of whom 3 are first-class, 1 is third-class and 3 are fourth-class. Suppose that 5 units of service
are available at time $n$. Then the departures at time $n$ will be 3 first-class customers,
1 third-class customer and 1 fourth-class customer, leaving 2 fourth-class customers
remaining in the queue. 

Let $m$ be the total number of classes. We will have a collection of arrival processes
$A^{=1}, A^{=2},\dots,A^{=m}$, where $A^{=r}(n)$ is the number of $r$th-class customers
arriving at time $n$. 

We will also denote $A^{\leq r}=\sum_{i=1}^r A^{=i}$, for $r=0,1,\dots,m$. 

Let $S$ be the service process of the queue. 

Similarly we will write $D^{=r}(n)$ for the number of $r$th-class customers departing
at time $n$, and $X^{=r}(n)$ for the number of $r$th-class customers present in the 
queue at the beginning of time-slot $n$. Write also
$D^{\leq r}=\sum_{i=1}^r D^{=i}$ and $X^{\leq r}=\sum_{i=1}^r X^{=i}$.

There are two natural ways to construct the multi-class queue (which are equivalent).
One way would be to look at the queueing process of $r$th-class customers for each $r$.
This is a queue with arrival process $A^{=r}$ and service process $S-D^{\leq r}$; 
the services available to $r$th-class customers are those that have not been 
used by any higher-priority customer. 

Alternatively, we will consider, for each $r$, the queueing process of customers
of classes $1,2,\dots, r$ combined. This is a queue with arrival process $A^{\leq r}$ and
service process $S$. So in particular we define
\begin{align*}
D^{\leq r}&=D(A^{\leq r}, S),\\
X^{\leq r}&=X(A^{\leq r}, S).
\end{align*}

This second description turns out to be more useful since it describes the $m$-class
queue as a coupling of $m$ single-class queues, each with the same service process. 
Once we come to consider interchangeability of queues, the fact that 
we are working with a common service process becomes crucial. 

\section{Bernoulli-geometric distribution}\label{BerGeom}\label{BerGeomqueue}
We define a \textit{Bernoulli-geometric distribution},
with parameters $p$ and $\alpha$.
A random variable with this distribution has the distribution of the product of 
two independent random variables, one with $\Ber(p)$ distribution and the other with 
$\Geom(\alpha)$ distribution. That is, 
$A\sim\Ber(p)\Geom(\alpha)$ if 
\[
\PP(A=k)=\begin{cases}
1-p,& k=0\\
p\alpha(1-\alpha)^{k-1},& k\geq 1.\end{cases}
\]
We have $\EE A=p/\alpha$.

We will consider a queue where $A(n)$ is an i.i.d.\ sequence with $A(n)\sim\Ber(p)\Geom(\alpha)$,
and $S(n)$ is an i.i.d.\ sequence (independent of $A(n)$) with $S(n)\sim\Ber(q)\Geom(\beta)$.
(We will say that $A$ and $S$ are \textit{Bernoulli-geometric processes}).
Queues of this type are investigated in \cite{jbmbatch}, 
in particular regarding their reversibility properties. 

For stability we will assume that $\EE A(n)<\EE S(n)$, i.e.\ that $p/\alpha<q/\beta$,
and further we assume that
\begin{equation}\label{1param}
c(p,\alpha)=c(q,\beta),
\end{equation}
where for $p,\alpha\in[0,1]$ we define 
\begin{equation}\label{cdef}
c(p,\alpha)=\frac{p}{1-p}\frac{\alpha}{1-\alpha}.
\end{equation}

When the arrival and service processes share a value of the parameter $c$ in this way, 
a version of Burke's theorem holds.
\cite{jbmbatch}. We use an asterisk to denote the time reversal of a process, 
so that $A^*(n)=A(-n)$. 
The following result is included in Theorem 4.1 of \cite{jbmbatch}:

\begin{theorem}\label{Burkethm}
Under the assumptions above, the departure process $\left(D(n), n\in\ZZ\right)$ has the same
law as the arrival process $\left(A(n), n\in\ZZ\right)$.
Moreover, the queue is reversible in the sense that 
\[
\left(A(n), D(n),n\in\ZZ\right) \isd \left(D^*(n), A^*(n), n\in\ZZ\right).
\]
\end{theorem}

Note that for given $q$, $\beta$ and $\lambda$, there is precisely one choice of 
the parameters $p$ and $\alpha$ such that $p/\alpha=\lambda$ and $c(p,\alpha)=c(q,\beta)$.
So for a Bernoulli-Geometric service process with given parameters $q$ and $\beta$, 
Theorem \ref{Burkethm} gives a one-parameter family of fixed-point arrival processes,
indexed by the arrival intensity $p/\alpha$. 

In the next section we see further that two service processes are interchangeable
precisely when they fall within the same one-parameter family (that is, 
when they share a value of the parameter $c$). 

\section{Interchangeability}\label{interchangeability}
We consider two queues in tandem. Suppose we are given an arrival process $A_1$ 
and service processes $S_1, S_2$.
From this we construct a system of two queues.
The first queue has arrival process $A_1$ and service process $S_1$. 
The second queue has service process $S_2$ and arrival process $A_2$ defined by 
$A_2=D(A_1, S_1)$. (The departure process from the first queue
becomes the arrival process from the second queue). 

Note that a departure from queue 1 at time $n$ becomes an arrival at queue 2
at the same time-slot $n$ (and may indeed depart from queue 2 at time $n$ also).

The following result says that two service processes, which fall within the
same one-parameter family mentioned above, are \textit{interchangeable}:

\begin{theorem}\label{interchangethm}
Suppose that the processes $S_1$ and $ S_2$ are independent of each other
and each are i.i.d., with $S_1(n)\sim\Ber(p)\Geom(\alpha)$ and $ S_2(n)\sim\Ber(q)\Geom(\beta)$.
Suppose further that condition (\ref{1param}) holds. 

Then there is a coupling of $(S_1,  S_2)$ with another pair of processes $(\tS_1, \tS_2)$ such that:
\begin{itemize}
\item[(i)] $(S_1, S_2)\isd(\tS_2, \tS_1)$;
\item[(ii)] If $A$ is any arrival process, then $D_2=\tD_2$ where
$D_2=D(D(A,S_1),  S_2)$ and $\tD_2=D(D(A,\tS_1), \tS_2)$. 
\end{itemize}
\end{theorem}

Note that in particular it follows that for any $A$, 
the distributions of $D(D(A,S_1),  S_2)$ and $D(D(A, S_2), S_1)$ are the same.
The result gives something rather stronger; namely a coupling of 
the two pairs of service processes such that, simultaneously for all arrival processes,
the output of the system is the same for both pairs. 

Recall that a multitype queue can be seen as a coupling of several queues,
with different arrival processes but the same service process. Because 
Theorem \ref{interchangethm} gives a coupling which works for all 
arrival processes simultaneously, we immediately have as a corollary an 
interchangeability result for multitype queues. 

\begin{theorem}\label{multitypeinterchangethm}
Suppose the conditions of Theorem \ref{interchangethm} hold, and
let $\uA$ be a multitype arrival process. 
Then $D(D(\uA, S_1), S_2)\isd D(D(\uA, S_2), S_1)$.
\end{theorem}

\section{Multitype fixed points}\label{mainresultsection}
In this section we will define a series of distributions of multi-class processes,
denoted by $F_{m,c,\lambda_1,\dots,\lambda_m}$.

$F_{m,c,\lambda_1,\dots,\lambda_m}$ will be the distribution of an $m$-type process,
with intensity of $j$th-class customers equal to $\lambda_j$. We will show that it is a fixed point
arrival process for a $\Ber(q,\beta)$ service process, for all $q$ and $\beta$ such
that $q/\beta\geq\lambda_1+\dots+\lambda_m$ and $c(q,\beta)=c$ (where the function
$c$ is defined at (\ref{cdef})).
When $c$ and $\{\lambda_j\}$ are fixed, we may abbreviate and write simply $F_m$.
The distributions are defined recursively, using $F_{m-1}$ to construct $F_m$.

$F_1$ is the distribution of a 1-type process; simply a $\Ber(p)\Geom(\alpha)$ 
process where $p$ and $\alpha$ are chosen so that $p/\alpha=\lambda_1$ 
and $c(p,\alpha)=c$.

Suppose we have constructed $F_{m-1}$. Now $F_m$ is constructed as follows.
We consider a queue whose service process is $\Ber(q_m)\Geom(\beta_m)$,
where $q_m$ and $\beta_m$ are chosen so that the total service intensity $q_m/\beta_m$ 
is equal to $\lambda_1+\dots+\lambda_m$, and so that $c(q_m,\beta_m)=c$.
(As observed in Section \ref{BerGeom}, there is a unique such choice).

The arrival process to the queue has distribution $F_{m-1}$ (and is independent
of the service process). This leads to an $(m-1)$-type departure process. 
We will extend the departure process to an $m$-type process by 
replacing the unused service in the queue by customers of type $m$. 
Let $F_m$ be the distribution of the $m$-type process obtained in this way.

Note that a consequence of the construction and of the fixed point result is 
that the restriction of a process with distribution $F_m$ to its first $m-1$ 
types gives a process with distribution $F_{m-1}$.

\begin{theorem}\label{fixedpointtheorem}
Let $S$ be a $\Ber(q)\Geom(\beta)$ service process, with $q/\beta=\mu$ and $c(q,\beta)=c$.
Then for all $\{\lambda_j\}$ with $\lambda_1+\dots+\lambda_m\leq \mu$,
the distribution $F_{m,c,\lambda_1,\dots,\lambda_m}$ is a fixed point
for the process $S$.

That is, suppose $\uA$ is an $m$-type arrival process with distribution 
$F_{m,c,\lambda_1,\dots,\lambda_m}$ and is independent of $S$. 
If $\uD=D(\uA, S)$ is the $m$-type departure process, then $\uD\isd\uA$.
\end{theorem}

\noindent\textit{Proof}:
The case $m=1$ is the version of Burke's theorem given in Theorem \ref{Burkethm} above.
The arrival process is a one-type $\Ber(p)\Geom(\alpha)$ process and
(\ref{1param}) holds, so indeed the departure process has the same law as the arrival process.

Now let $m\geq2$ and suppose that the result holds as claimed for $F_{m-1}$. We wish to establish
the result for $F_m$. We will use the interchangeability result 
in Theorem \ref{multitypeinterchangethm}.

We consider two independent service processes $S_\mu$ and $S_{\lambda_1+\dots+\lambda_m}$, 
whch are Bernoulli-geometric processes with parameters $\Ber(r,\gamma)$ and $\Ber(q_m,\beta_m)$
respectively. We choose the parameters such that:
\begin{itemize}
\item $r/\gamma=\mu$ and $c(r,\gamma)=c$;
\item $q_m/\beta_m=\lambda_1+\dots+\lambda_m$ and $c(q_m,\beta_m)=c$.
\end{itemize}

We will consider an $m$-type arrival process $\uA$ distributed according to $F_m$.
We consider the effect of feeding this arrival process into two tandem systems, 
comprising the service processes $S_\mu$ and $S_{\lambda_1+\dots+\lambda_m}$ in their two possible orders.
Theorem \ref{multitypeinterchangethm} tells us that $D(D(\uA, S_\mu), S_{\lambda_1+\dots+\lambda_m})$ and 
$D(D(\uA, S_{\lambda_1+\dots+\lambda_m}), S_\mu)$ have the same distribution.
To make use of this result we first need two simple properties. 
The first concerns projections of multi-type processes onto processes
with fewer types, and the second concerns the behaviour of a queue which is saturated
(i.e.\ whose service rate is only just sufficient to serve all the arriving customers).

\textit{Claim 1:}
Suppose $\uA$ has distribution $F_m$, and let $G_m$ be the distribution of $D(\uA, S_{\mu})$.
Then the restrictions of $F_m$ and $G_m$ to the first $m-1$ coordinates 
have the same distribution, and the combined process of all $m$ types
under $G_m$ has the same distribution as $S_{\lambda_1+\dots+\lambda_m}$.

\textit{Proof of Claim 1:}
Restricting $F_m$ to its first $m-1$ coordinates gives
the distribution $F_{m-1}$. Since $F_{m-1}$ is an $(m-1)$-type fixed
point for the service process $S_{\mu}$, the first $m-1$ types
under $G_m$ do indeed have distribution $F_{m-1}$ also. 

Finally, from the defintion of $F_m$, the combined process of all $m$ types under $F_m$ is
a $\Ber(q_m)\Geom(\beta_m)$ process.   
The service process $S_{\mu}$ is a $\Ber(r)\Geom(\gamma)$ process. 
By assumption $c(q_m, \beta_m)=c(r,\gamma)=c$. So 
condition (\ref{1param}) is satisfied for the queue, and we can apply 
the 1-type fixed point result in Theorem \ref{Burkethm} to show that 
the distribution of the combined process of all $m$ types in the departure process
has the same distribution as the combined process of all $m$ types in the arrival
process, namely $S_{\lambda_1+\dots+\lambda_m}$. This completes the argument for Claim 1.

\textit{Claim 2:}
Suppose $\uA$, restricted to its first $m-1$ coordinates, has
distribution $F_{m-1}$. Suppose also that the combined process of all $m$
types of arrival in $\uA$ is an i.i.d.\ process with intensity 
$\lambda_1+\dots+\lambda_m$. 
Then $D(\uA, S_{\lambda_1+\dots+\lambda_m})$ has
distribution $F_m$. 

\textit{Proof of Claim 2:}
Recall that $F_m$ is the distribution obtained by passing an
$(m-1)$-type arrival process with distribution $F_{m-1}$ through a
queue with service process $S_{\lambda_1+\dots+\lambda_m}$, 
and putting customers of type $m$ in place of all unused service. 

The first $m-1$ components of $\uA$ indeed have distribution $F_{m-1}$. 
Certainly the $m$th component of the departure process is
a subset of the unused service from the $m-1$ first types. 
To show that the departure process has distribution $F_m$, it remains
to show that all the service unused by the first $m-1$ types 
is used by customers of type $m$; that is, that if we look
at all $m$ types combined, then there is no unused service.

The arrival process of all customers combined is an i.i.d.\ process
with rate $\lambda_1+\dots+\lambda_m$; the same is true of the 
service process, and the arrival and service processes are independent. 
Hence by Lemma \ref{saturationlemma}, there is no unused service
in the queue as desired. So indeed the output process has distribution $F_m$,
and we have established Claim 2 as required.

Now we put together Claims 1 and 2 to complete the proof. We consider
the two possible orderings of the service processes $S_\mu$ and 
$S_{\lambda_1+\dots+\lambda_m}$. The two systems that arise can be illustrated as follows:
\begin{align*}
F_m\longrightarrow\boxed{\hspace{6.09mm}S_\mu\hspace{6.09mm}}
&\longrightarrow \makebox[1.5em]{$G_m$} \longrightarrow
\boxed{S_{\lambda_1+\dots+\lambda_m}}
\longrightarrow{H_m}\\
F_m\longrightarrow\boxed{S_{\lambda_1+\dots+\lambda_m}}
&\longrightarrow \makebox[1.5em]{$J_m$} \longrightarrow
\boxed{\hspace{6.09mm}S_{\mu}\hspace{6.09mm}}
\longrightarrow{K_m}
\end{align*}

By Claim 1, $G_m$ satisfies the conclusion of Claim 1, which is also
the assumption on the arrival process of Claim 2. Then by Claim 2, $H_m=F_m$.

Certainly $F_m$ itself also satisfies the assumption of Claim 2,
so by Claim 2, $J_m=F_m$. 

Theorem \ref{interchangethm} says that $H_m=K_m$. 

Then the second box of the second line down shows us that
if we pass $J_m$ through $S_\mu$ we get $K_m$.

But $J_m=K_m=F_m$. So indeed $F_m$ is a fixed point for $S_\mu$ as desired.$\hfill\Box$

\section{Proof of interchangeability result}\label{interchangeabilityproof}
In this section we will prove Theorem \ref{interchangethm}. 
We first need two simple lemmas about queues in tandem.

\begin{lemma}\label{trunclemma}
Consider a system of two queues in tandem with arrival process $A$ and service processes
$S_1$ and $S_2$, all taking integer values. For $s\in\ZZ$, define $A|_{(-s, \infty)}$ to be the 
arrival process truncated before time $-s$. That is, 
\[
A|_{(-s, \infty)}(n)=\begin{cases} 0&\text{ if $n< -s$},\\A(n)&\text{ if $n\geq -s$.}\end{cases}
\]
Define $D_2=D(D(A,S_1), S_2)$ and $D_2^{(s)}=D(D(A|_{(-s,\infty)}, S_1), S_2)$. 
For any $n$, if $s$ is large enough then $D_2(n)=D_2^{(s)}(n)$.
\end{lemma}
\noindent\textit{Proof:}
We compare the original system and the system with the arrival process
before 
time $-s$. Write $D_1$ and $D_2$ for the departure processes from the first and second queues 
in the original system, and similarly $D_1^{(s)}$, $D_2^{(s)}$ for 
the processes in the system with the truncated arrival process.

Take any $n\in\ZZ$.
From the definitions in Section \ref{model}, we have
\begin{align*}
D_1(n)&=\min\left(
S_1(n),A(n)+\sup_{u\leq n}\sum_{r=u}^{n-1}\left[A(r)-S_1(r)\right]\right),\\
D_1^{(s)}(n)&=\min\left(
S_1(n),A(n)+\sup_{-s\leq u\leq n}\sum_{r=u}^{n-1}\left[A(r)-S_1(r)\right]\right).
\end{align*}
Certainly $D_1^{(s)}(n)$ is increasing in $s$, and $D_1^{(s)}(n)\leq D_1(n)$ for all $s$.

We first observe that $D_1^{(s)}(n)=D_1(n)$ for all large enough $s$. 
Suppose the sup in the first line is finite. Then (since the variables all take 
integer values) it is attained for some $u=u^*$, and we have $D_1(n)=D^{(s)}_1(n)$ for 
all $s>-u^*$. If instead the sup in the first line is infinite, then 
by taking $s$ large we may make the sup in the second line as large as desired; 
if $s$ is large enough that the second term in the min exceeds $S_1(n)$, 
then we have $D_1(n)=D^{(s)}_1(n)$, as required.

Now we have the corresponding expressions for the departure process from the second queue
(whose arrival process is $D_1$ and whose service process is $S_2$):
\begin{align*}
D_2(n)&=\min\left(
S_2(n),D_1(n)+\sup_{u\leq n}\sum_{r=u}^{n-1}\left[D_1(r)-S_2(r)\right]\right),\\
D_2^{(s)}(n)&=\min\left(
S_2(n),D_1^{(s)}(n)+\sup_{u\leq n}\sum_{r=u}^{n-1}\left[D_1^{(s)}(r)-S_2(r)\right]\right).
\end{align*}
Since 
$D_1^{(s)}(n)$ is increasing in $s$ and bounded above by $D_1(n)$, we also  
have that $D_2^{(s)}(n)$ is increasing in $s$ and bounded above by $D_2(n)$. 

Suppose that the sup in the first line is finite. Then as before
it is attained at some $u^*$. Now if we take $s$ large enough that
$D_1^{(s)}(r)=D_1(r)$ for all $r$ with 
$u^*\leq r\leq n$, then we see from the second line
that indeed $D_2^{(s)}(n)\geq D_2(n)$, and hence in fact $D_2^{(s)}(n)=D_2(n)$ as required.
A similar argument applies if the sup in the first line is infinite;
by taking $s$ large enough, we can make the sup in the second line as large as desired.$\hfill\Box$

\begin{lemma}\label{queuereplemma}
Consider a tandem of two queues, with arrival process $A$ and service processes $S_1$ and $S_2$.
Suppose that $A(n)=0$ for all $n<0$. Let $D_2$ be the departure process from the second queue.
Then for all $t\geq 1$,
\begin{equation}\label{queuerepformula}
\sum_{r=0}^{t-1} D(r)
=\inf_{0\leq u_1\leq u_2\leq t}
\left\{
\sum_{r=0}^{u_1-1} A(r)
+\sum_{r=u_1}^{u_2-1} S_1(r)
+\sum_{r=u_2}^{t-1} S_2(r)
\right\}
.
\end{equation}
\end{lemma}
\noindent\textit{Proof:}
An equivalent result for a related queueing model was given in \cite{SzcKel}, 
and similar properties appear in many places in the literature. 
Formula \ref{queuerepformula} is essentially a special case of equation (10) of 
\cite{OConnell}.

Let $D_1$ and $D_2$ be the departure processes from the first and second queues,
and $X_1$ and $X_2$ be the queue-length processes at the first and second queues.

Since $A_n=0$ for all negative $n$, the queues start empty at time 0: $X_1(0)=X_2(0)=0$. 
So the total number of departures from the system (that is, from the second queue)
before time $t$ is given by the total number of arrivals during that time, minus
the queue-lengths at time $t$:
\begin{equation}\label{Dsum}
\sum_{r=0}^{t-1}D_2(r)=\sum_{r=0}^{t-1}A(r) - X_1(t) - X_2(t).
\end{equation}

Similarly, considering the evolution of the first queue alone between times $u$ and $t$, 
we have that for any $u<t$,
\[
\sum_{r=u}^{t-1}D_1(r)=\sum_{r=u}^{t-1}A(r) - X_1(t) + X_1(u).
\]

Now $X_2$ is the queue-length of a queue with arrivals $D_1$ and services $S_2$. 
So from the definition (\ref{Xdef}), we have
\begin{align*}
X_1(t)+X_2(t)
&
=X_1(t)+\sup_{u_2\leq t} \sum_{r=u_2}^{t-1} \left[D_1(r)-S_2(r)\right]\\
&
=X_1(t)+\sup_{u_2\leq t}
\left\{
\sum_{r=u_2}^{t-1}A(r) 
- X_1(t) + X_1(u) 
-\sum_{r=u_2}^{t-1}S_2(r)\right\}\\
&=\sup_{u_2\leq t}
\left\{
X_1(u)+\sum_{r=u_2}^{t-1}
\left
[A(r)-S_2(r)
\right]
\right\}\\
&
=\sup_{u_2\leq t}
\left\{
\sup_{u_1\leq u_2} \sum_{r=u_1}^{u_2-1}
\left[A(r)-S_1(r)\right]+
\sum_{r=u_2}^{t-1}\left[A(r)-S_2(r)\right]\right\}\\
&=\sup_{u_1\leq u_2\leq t}\left\{
\sum_{r=u_1}^{t-1}A(r) 
-\sum_{r=u_1}^{u_2-1} S_1(r)
-\sum_{r=u_2}^{t-1} S_2(r)\right\}.
\end{align*}
Since $A_n=0$ for negative $n$, the sup will be attained for some $u_1\geq 0$. 
Then using (\ref{Dsum}) gives (\ref{queuerepformula}) as desired.$\hfill\Box$

Now we proceed to prove Theorem \ref{interchangethm}.
To do this, 
we need to construct $S_1,  S_2, \tS_1$ and $\tS_2$ on the same probability space, 
in such a way that:
\begin{itemize}
\item[(i)]$(S_1, S_2)\isd(\tS_2,\tS_1)$
\item[(ii)] For all $A$, $D(D(A,S_1), S_2)=D(D(A,\tS_1),\tS_2)$.
\end{itemize}

We assume that $S_2$ has higher intensity than $S_1$. (If not, then swap them around. If the intensity
is the same, then the assumptions ensure that they have the same distribution.)

The construction is as follows. We think of a queue with arrival process $S_1$ and service
process $S_2$. We can decompose $S_2$ into $D(S_1,S_2)$ and $U(S_1,S_2)$ -- that is, 
$S_2=D(S_1,S_2)+U(S_1,S_2)$. 

We let $\tS_2=D(S_1,S_2)$. So we have obtained $\tS_2$ from $ S_2$ by removing the unused
service process. By Theorem \ref{Burkethm}, 
we have
that $(S_1, \tS_2)\isd(\tS_2^*,S_1^*)$. (Recall that the asterisk denotes the reverse of a process).
Knowing this, we can extend the probability space to include a random variable $\tS_1$, with 
its time-reversal $\tS_1^*$,  such that
$(S_1, \tS_2,  S_2)\isd(\tS_2^*,S_1^*,\tS_1^*)$. 

In particular this implies that $(S_1, S_2)\isd(\tS_2^*, \tS_1^*)$. But 
the distribution of $(S_1, S_2)$ is invariant under time reversal, so this definition
gives us (i) above as desired.

(Note that since $(S_1, \tS_2,  S_2)\isd(\tS_2^*,S_1^*,\tS_1^*)$, we have in particular that
$S_1^*=D(\tS_2^*, \tS_1^*)$, 
or equivalently that $S_1^*+U(\tS_2^*,\tS_1^*)=\tS_1^*$.
So we can think of $\tS_1$ as being obtained from $S_1$ by 
adding an unused service process, for a queue operating in reverse time). 

Now we need to verify (ii) above. Using Lemma \ref{trunclemma}, 
it will be sufficient to verify 
it for every arrival process which is 0 up to some finite time. 
Without loss of generality, suppose that $A(n)=0$ for all $n<0$.
Let $D=D(D(A,S_1),S_2))$.
From Lemma \ref{queuereplemma}, we have the 
following representation of the total number of departures
from the tandem of two queues from time 0 up to time $t-1$:
\[
\sum_{r=0}^{t-1} D(r)
=\inf_{0\leq u_1\leq u_2\leq t}
\left\{
\sum_{r=0}^{u_1-1} A(r)
+\sum_{r=u_1}^{u_2-1} S_1(r)
+\sum_{r=u_2}^{t-1} S_2(r)
\right\}
.
\]
It will be enough to show that for any $A$, this quantity is unchanged if we replace
$(S_1,S_2)$ by $(\tS_1, \tS_2)$. 
For this, in turn it's enough that for all $s$ and $t$, the quantity
\[
\inf_{s\leq u\leq t}
\left\{
\sum_{r=s}^{u-1} S_1(r) + \sum_{r=u}^{t-1} S_2(r)
\right\}
\]
is unchanged if we replace $(S_1,S_2)$ by $(\tS_1,\tS_2)$. 

We consider the queue with arrival process $S_1$ and service process $S_2$,
along with its queue-length, unused service and departure processes $X$, $U$ and $D$
which are functions of $S_1$ and $S_2$. Note that $U(n)>0$ precisely if $\tS_2(n)<S_2(n)$.

First we show that the quantity we are interested in 
is unchanged if we replace $(S_1,S_2)$ by $(S_1, \tS_2)$; that is:
\begin{equation}\label{compare}
\inf_{s\leq u\leq t}
\left\{
\sum_{r=s}^{u-1} S_1(r) + \sum_{r=u}^{t-1} S_2(r)
\right\}
=
\inf_{s\leq u\leq t}
\left\{
\sum_{r=s}^{u-1} S_1(r) + \sum_{r=u}^{t-1} \tS_2(r)
\right\}.
\end{equation}

If there is no unused service between times 0 and $t-1$, 
then $S_2$ and $\tS_2$ agree on the whole interval and the equality is obvious.

Otherwise, let $n$ be the latest time in $\{s,\dots,t-1\}$ for which there is unused service,
i.e.\ for which $S_2(n)>\tS_2(n)$, or equivalently that $U(n)>0$ where $U=U(S_1,S_2)$. 
In particular there are no more arrivals than departures in the interval $\{s,\dots,n\}$
(otherwise there would still be customers in the queue after time $n$, and hence
there could not have been unused service at time $n$).
That is, 
$\sum_{i=s}^n S_1(i)\leq \sum_{i=s}^n \tS_2(i)$. 
This implies that we can take $u>n$ in minimising the RHS of (\ref{compare}). 
Then consider using the same $u$ on the LHS as on the RHS. Since $u>n$ and
$n$ was the last moment of unused service, the processes $S_2$ and $\tS_2$ agree 
after time $n$. So in fact the LHS is at least as small as the RHS. 
But since $S_2\geq \tS_2$ it is also clear that the LHS is no smaller than the RHS. 
So in fact the two are the same. 

Now we want to show that we can further replace $S_1$ by $\tS_1$; that is,
\begin{equation}
\inf_{s\leq u\leq t}
\left\{
\sum_{r=s}^{u-1} S_1(r) + \sum_{r=u}^{t-1} \tS_2(r)
\right\}
=
\inf_{s\leq u\leq t}
\left\{
\sum_{r=s}^{u-1} \tS_1(r) + \sum_{r=u}^{t-1} \tS_2(r)
\right\}.
\end{equation}
This is in fact precisely the same relation as (\ref{compare}), but now for the queue
with arrival process $\tS_2^*$ and service process $\tS_1^*$. So an exactly analogous
argument can be applied. 

This completes the proof of (ii) and hence of Theorem \ref{interchangethm}.

\section{Summary and examples}\label{examplesection}
In this section we point out various interesting cases that can be 
obtained by taking particular values or appropriate limits in 
the general model considered up to now. 

Our general model is of a service process in discrete time, whose batches
are i.i.d.\ with $\Ber(q)\Geom(\beta)$ distribution. Such a process
has service intensity $q/\beta$.
This model has two 
parameters. For each choice $(q,\beta)$ there is a one-dimensional 
family of service process which are interchangeable with it; namely
an $(r,\gamma)$-process is interchangeable with a $(q,\beta)$-process
whenever $c(r,\gamma)=c(q,\beta)$. This relation partitions the parameter space
into equivalence classes (and each class contains precisely one 
set of parameters for each service intensity).

Suppose we are given service parameters $(q,\beta)$ with service intensity
$\mu=q/\beta$, and also $\lambda_1,\dots,\lambda_m$ such that
$\lambda_1+\dots+\lambda_m<\mu$. Then there is precisely
one $m$-type fixed point arrival process with intensities $\lambda_r$ of $r$th-class
customers for $1\leq r\leq m$. For example, the two-type fixed point with 
intensities $\lambda_1$ and $\lambda_2$ is obtained by taking
the departure and unused service processes from a queue with arrival process $A$ 
and service process $S$ where $A$ and $S$ are chosen so that $A$ has intensity 
$\lambda_1$, $S$ has intensity $\lambda_1+\lambda_2$, and both $A$ and $S$
are interchangeable with the $(q,\beta)$-process. (As explained in the last
paragraph, this uniquely determines $A$ and $S$). 

This framework extends to various particular cases which we describe below. 
In some cases, the model has only one parameter, and any pair of 
processes in the class considered are interchangeable. 

\subsection{$./M/1$ queue in discrete time}\label{discretedotM1section}
If we take $\beta=1$, the service process is a ``Bernoulli process''. 
The batches are i.i.d., each equal to 0 with probability $1-q$ and 1 
with probability $q$, so that $\mu=q$.  
Hence we have a one-parameter model. 

This is the discrete-time equivalent of the $./M/1$ server described in 
Section \ref{dotM1section}. A version of Burke's theorem for such systems
was proved by Hsu and Burke in \cite{HsuBurke}; any Bernoulli process with
lower intensity $p<q$ is a one-type fixed point arrival process. 

Any two such Bernoulli processes are interchangeable. 
For any $m$ and $\lambda_1+\dots+\lambda_m<1$, there
is an arrival process $F_{m,\lambda_1,\dots,\lambda_m}$ 
with intensities $\lambda_i$ and which is a fixed
for any such $./M(\mu)/1$ server whenever $\mu>\lambda_1+\dots+\lambda_m$. 
The multi-type fixed points
were constructed in \cite{FerMarmulti}, and they can also be seen
as stationary distributions for multi-type versions of the TASEP
(totally asymmetric simple exclusion process). 

\subsection{Geometric and exponential batches in discrete time}
If we take $q=1$, then the service batches are i.i.d.\ geometric
distributions. (A similar case arises for $q=1-\beta$, 
except that the geometric distribution obtained starts from 0 rather than from 1).
Again we have a one-parameter model, and
any two such geometric processes are interchangeable. 
A related result (for a particular form of the arrival process)
was given by Draief, O'Connell and Mairesse in \cite{DraMaiOco}.
A version of Burke's theorem had been given by Bedekar and Azizo{\~g}lu
in \cite{BedAzi}.

The multi-type fixed points in this case are related to stationary 
distributions for a multi-type version of the totally asymmetric zero-range process
\cite{jbmTAZRP}.

Now take the limit $\delta\to0$, putting $\beta=\delta/u$ 
and rescaling work by a factor $\delta$. 
In this way 
we can obtain the case where the batches are i.i.d.\ exponential 
with mean $u$. Again, any two such exponential processes are interchangeable. 

\subsection{Batch queues in continuous time}
Suppose we let $\epsilon\to0$, take 
$q=\nu\epsilon$ and rescale time by a factor $\epsilon$. 

Then we have batches which are geometric with parameter
$\beta$, which occur at times of a Poisson process of rate $\nu$.

Now we have a two-parameter distribution, and for each pair of parameters
$(\nu, \beta)$ there is a one-parameter family of processes, one
for each service intensity, which are interchangeable with it. 
Namely, processes of this kind with parameters $(\nu,\beta)$ and $(\rho,\gamma)$
are interchangeable if 
\[
\frac{\nu\beta}{1-\beta}=\frac{\rho\gamma}{1-\gamma}.
\]

If we further take $\delta\to0$, putting $\beta=\delta/u$ and 
rescaling work by a factor $\delta$, then
we obtain the case where the batches are exponential with mean $u$ (again
occurring at the times of a Poisson process of rate $\nu$). Two
such ``Poisson-exponential'' processes with parameters $(\nu, u)$ and $(\rho, v)$ are 
interchangeable if $\nu/u$=$\rho/v$.

If on the other hand we take $\beta\to 1$, the batches all have size 1 and 
we obtain the familiar $./M/1$ queueing server which we described in detail in 
Section \ref{dotM1section}.

\subsection{Brownian queues}
Queues in continuous time (as in the previous example) can be represented using a
notation analogous to that introduced for discrete-time queues in Section \ref{model}.

Let $S_t,t\in\RR$ be a service process and $A_t,t\in\RR$ be an arrival process.
We interpret $S_t-S_s$ as the amount of service offered in the interval $(s,t]$,
and similarly $A_t-A_s$ as the amount of work arriving in $(s,t]$.

Adding a constant to either of the processes $A$ and $S$ makes no difference,
so we may take for example $A_0=S_0=0$.

We define processes $Q_t$, $D_t$ and $U_t$, as functions of the processes $A_t$ and $S_t$, by
\begin{align*}
Q_t&=\sup_{s<t}\Big\{ \big(A_t-A_s)-\big(S_t-S_s) \Big\},\\
D_t&=A_t+Q_0-Q_t,\\
U_t&=S_t-D_t.
\end{align*}

Now $D_t-D_s$ is the amount of work departing in $(s,t]$, $U_t-U_s$ is the amount
of unused service in $(s,t]$, and $Q_t$ is the queue-length at time $t$.

Note that under these definitions we have $U_0=D_0=0$ (again, this normalization is not important).

A frequently studied example is that of the \textit{Brownian queue}. Let $A_t$ be a two-sided
Brownian motion with drift $\lambda$ and variance 1, and let $S_t$ be a two-sided
Brownian motion with drift $\mu$ and variance 1, where $\lambda, \mu\in\RR$ and $\lambda<\mu$. 

This model arises naturally as the scaling limit of queues in the so-called ``heavy traffic'' regime
\cite{Harrisonbook}, \cite{Williamssurvey}.
For example, for $n\in\NN$ consider a discrete time queue whose arrival process $A(r),r\in\ZZ$
is a Bernoulli process with rate $1/2+\lambda/\sqrt{n}$ and whose service process $S(r), r\in\ZZ$
is a Bernoulli process with rate $1/2+\mu/\sqrt{n}$. We recentre and rescale by defining, for $t\in\RR$,
\begin{align*}
\tilde{A}^{(n)}_t&=\frac{\sum_{r=0}^{nt}A(r)-\frac{nt}{2}}{\sqrt{n/2}},
\tilde{S}^{(n)}_t&=\frac{\sum_{r=0}^{nt}S(r)-\frac{nt}{2}}{\sqrt{n/2}}.
\end{align*}
Then as $n\to\infty$, the distribution of the processes $\tilde{A}^{(n)}_t$ and $\tilde{S}^{(n)}_t$
converges to that of $A_t$ and $S_t$ above. 

A version of Burke's theorem holds for such a Brownian queue: namely that the processes $D_t$ 
and $A_t$ have the same distribution. 
This was shown by Harrison and Williams in 
\cite{HarWilBurke} (see also \cite{OcoYor} in particular, where the result is generalized in a 
number of ways). The process $Q_t$ has the law of a stationary reflected Brownian motion
with drift $-(\mu-\lambda)$, and its stationary distribution is exponential with rate $\mu-\lambda$.
The process $U_t$ is the local-time process of $Q_t$ at 0. Note that $U$ is non-decreasing,
and grows at average rate $\mu-\lambda$, but is constant except on a set of measure $0$
(namely, the set of times when $Q=0$). 

Thus we may say that a Brownian queueing server (with drift $\mu$) has Brownian motions of drift
$\lambda<\mu$ as fixed point arrival processes. Again this can be extended to multi-type arrival 
processes.

Given a service process $S_t$ and a two-type arrival process $A^{(1)}_t, A^{(2)}_t$, 
the two-type departure process $D^{(1)}_t, D^{(2)}_t$ is defined in a natural way by 
\begin{gather*}
D^{(1)}=D\left(A^{(1)}, S\right)
\\
D^{(1)}+D^{(2)}=D\left(A^{(1)}+A^{(2)},S\right).
\end{gather*}

Suppose we wish to construct a two-type fixed point for a Brownian service process of drift $\mu$,
which has drifts $\lambda_1$ and $\lambda_2$ where $\lambda_1<\lambda_1+\lambda_2<\mu$. 
To do this, consider a Brownian queue whose arrival process is a Brownian motion 
with drift $\lambda_1$, and whose service process is a Brownian motion with drift $\lambda_2$.
The process $(D,U)$ of departures and unused service from this queue is a two-type
fixed point of the required type. 

In this example we can see an extreme form of the ``clustering of lower-class customers''
already observed in the discrete case. Here the second-class work, represented
by the local-time process $U$, is concentrated on a set of times of measure 0. 
Note also the difference in nature between the process $D$ of first-class work
(which is non-monotonic, and in fact has unbounded variation) 
and the process $U$ of second-class work (which is non-decreasing). 

Fixed points with larger numbers of classes can be constructed recursively in the 
same way as in the discrete case. The subprocesses of $m$th-class work in these
fixed points are non-decreasing and singular, for each $m\geq 2$.

One can also obtain an interchangeability result; two
Brownian service processes, with different drifts but with the same variance, 
are interchangeable.

Arrival processes obtained as the local time of a reflecting diffusion 
have been considered before in various contexts, for example in 
\cite{ManNorSal} and \cite{KKSS}. These models were also motivated
in part by the modelling of a second-class departure process
as a local-time process. However, the models analysed are somewhat different; 
the service processes are deterministic rather than Brownian ($S_t=ct$)
and the local-time arrival process has access to all the service capacity
(rather than sharing the queue with a higher-priority stream).

\section{Interacting particle systems}\label{particlesection}
There is a close analogy between fixed-point arrival processes for queues and 
stationary distributions for interacting particle systems. 

Consider for example the TASEP (totally asymmetric simple exclusion process). In the 
one-type version of the process, each site of $\ZZ$ contains either a particle or a hole.
The state-space of the process can be written as $\{1,\infty\}^\ZZ$, where $1$ 
denotes a particle and $\infty$ denotes a hole. The dynamics of the process are as 
follows: each particle tries to jump to its left as a Poisson process of rate 1, 
and a jump succeeds if the site to the left of the particle is empty (in which case
the values $1$ and $\infty$ are exchanged between the two sites). 

One can also consider a multi-type TASEP. The process with $m$ types has
state space $(\{1,2\dots,m\}\cup\{\infty\})^\ZZ$. As before, each particle
tries to jump to its left as a Poisson process of rate 1, and a jump by 
a particle of type $r$ succeeds if the site to its left is occupied by a particle
with lower priority (higher-numbered class) or is empty (in which case 
the values at the two sites are interchanged as before). As in the case
of multi-type queues, this process can be seen in a natural way 
as a coupling of $m$ single-type TASEPs. 

Stationary distributions for TASEPs with two types of particle were constructed
by Derrida, Janowsky, Lebowitz and Speer \cite{DJLStasep}, and
related construction were also given in \cite{FFKtasep},  \cite{Angeltasep}
and \cite{DucSch}. In \cite{FerMarmulti}, this construction was extended
to multi-type processes, and made particularly explicit using a construction
using queues in tandem which corresponds to the construction
of multi-type fixed points for queues above.

Given $\lambda_1,\dots,\lambda_m$ with $\lambda_1+\dots+\lambda_m<1$, 
consider the fixed-point arrival process distribution $F_{m,\lambda_1,\dots,\lambda_m}$ 
for the discrete-time $./M/1$ queue
as described in Section \ref{discretedotM1section}. 
Now regard this distribution as a distribution over configurations of the $m$-type
TASEP; an arrival of type $r$ at time $n$ corresponds to a particle of type $r$ at site
$n$, and an empty arrival slot at time $n$ corresponds to a hole at site $n$. 

The result of \cite{FerMarmulti} is that $F_{m,\lambda_1,\dots,\lambda_m}$ 
is a stationary distribution for the $m$-type TASEP (indeed, 
it is the unique ergodic stationary distribution which has intensities $\lambda_r$).

Note that \textit{time} in the context of the queue corresponds to \textit{space}
in the context of the particle system. 

The analogy between queues and particle systems can be made even clearer
by considering particle systems evolving in discrete time rather than continuous time.
The same distributions $F_m$ are also stationary for certain versions
of the multi-type TASEP in discrete time \cite{MarSchdiscrete}. (Here many particles
may try to jump at the same time-step, and there are various possible natural
conventions for the order in which the jump attempts are to be processed). 
Now a queueing server, which uses a service process to transform an
arrival process into a departure process, can be seen as analogous to
a set of jump attempts at particular sites at a given time-step of the particle system, 
which transforms the particle configuration before that time-step
into the new particle configuration after that time-step. Successive updates
in the particle system correspond to consecutive queueing servers in a system 
of queues in tandem.

The fixed-point property for the distributions $F_m$ can be derived from 
the result for the particle system (as was done in \cite{FerMarmulti}). 
To see this, note that if we take an $m$-class TASEP and 
treat particles of class $m$ as holes,
then the resulting restriction of the process to the first $m-1$ types
is an $(m-1)$-class TASEP. Hence restricting the stationary distribution $F_m$
of the $m$-class system to its first $m-1$ types must give a stationary distribution
of the $(m-1)$-class system, namely $F_{m-1}$. 
That is, $F_{m-1}$ can be obtained from $F_m$ by restricting
to the first $m-1$ types. But also, the first $m-1$ types of $F_m$ are \textit{defined}
as the departure process from a queue whose arrival process is $F_{m-1}$ 
(the remaining $m$th class is the unused service process). So indeed $F_{m-1}$ 
is also the distribution of the departure process from the queue; this is the 
fixed-point result we obtained earlier.

However, such a proof via the particle system is rather indirect and perhaps
less natural than the approach using interchangeability above. 

In fact, in some cases, results previously obtained for stationary distributions
of particle systems can themselves be derived from the fixed-point results for queues.
For example, the stationary distributions for ``Hammersley's process'' \cite{AD95}
correspond to the fixed points for the $./M.1$ queue in continuous time,
described in Section \ref{dotM1section} \cite{FerMarHAD, femaihp}. 
The dynamics of Hammersley's 
process can be obtained by taking appropriate limits starting from
a $M/M/1$ queueing server in discrete time (first one lets the rates
of arrival and service approach 1 together in an appropriate way;
after rescaling space, the particles in a discrete-time version Hammersley's process
correspond to the gaps in the customer processes in the queue. 
Then the more familiar continuous-time version of Hammersley's proces
can be obtained by letting the rate of jumps in the discrete version go to 0 and 
rescaling time). In this way the result on stationary distributions for 
Hammersley's process can be seen as a consequence of the fixed-point results given here.
(However, we don't know how to use similar methods to derive the results
for the TASEP, for example). 

\section{Continuous class-labels}\label{continuoussection}
By considering the limit as the number of classes goes to infinity, and the 
intensity of each class goes to zero, one can arrive at models where
the classes are indexed by real numbers. 

This could be done in the general case above, but for simplicity we discuss
it in the case of a $./M/1$ server (where the service process is a Bernoulli process).
We will consider an arrival process in which at most one customer arrives 
at each time slot, and a customer arriving at time $n$ has a real-valued label. 
As before, lower-labelled customers have priority over higher-labelled customers. 

By taking appropriate limits in the results above, one can obtain the following result:

\begin{proposition}\label{continuousprop}
There exists a unique distribution over sequences $\big(L(n), n\in\ZZ\big)\in[0,1]^\ZZ$
with the following properties:
\begin{itemize}
\item[(i)]The distribution is stationary and ergodic.
\item[(ii)] $L(0)\sim U[0,1]$. 
\item[(iii)]Let $0<\lambda<1$. Define an arrival process as follows: if $L(n)\leq\lambda$,
then a customer with label $L(n)$ arrives at time $n$, while if $L(n)>\lambda$, 
then no customer arrives at time $n$. This arrival process is 
a fixed point for the $./M(\mu)/1$ queue for all $\mu>\lambda$.
\end{itemize}
\end{proposition}
Property (ii) is included just as a normalization (since the mechanics of the queue
are unchanged if the labels of all the customers are transformed by an increasing function). 
Then in property (iii), the condition $\mu>\lambda$ ensures that the queue is not saturated.

The distribution of $(L(n), n\in\ZZ)$ is also a stationary distribution 
for the TASEP. See \cite{AmiAngVal} for an investigation of this process, 
including in particular an interpretation as the ``speed process'' for a multi-type 
TASEP started out of equilibrium. 

One interesting property of the process is a manifestation of the ``clustering''
effect described above. Although $L(n)$ has a continuous distribution for each $n$,
nonetheless one has that for any $n$, $\PP(L(0)=L(n))>0$. 
(For example, $\PP(L(0)=L(1))=1/6$). In fact, with probability one there
exist infinitely many $n$ such that $L(0)=L(n)$. Hence clustering occurs in the following sense:
although any class-label has probability 0 of being seen a priori, if 
one sees the label at any particular time the same label has high probability of being seen
nearby, and will be seen infinitely often in the process.  

\section*{Acknowledgments}
JBM thanks Pablo Ferrari for many valuable conversations related to this work,
and Mike Harrison and Ilkka Norros for discussions about the results on Brownian queues.

\medskip

\parbox{0.33\textwidth}{\noindent
\textsc{Department of Statistics,\\
University of Oxford,\\
1 South Parks Road,\\
Oxford OX1 3TG,\\
UK}\\
\texttt{martin@stats.ox.ac.uk}\\$ $}
\hfill
\parbox{0.53\textwidth}{\noindent
\textsc{Departments of Electrical Engineering and Computer Science,\\
Packard Building,\\
Stanford University,\\
Stanford, CA 94305-9510,\\
USA}\\
\texttt{balaji@stanford.edu}}

\end{document}